\documentclass[12pt]{article}
\usepackage{mathrsfs}
\usepackage{amssymb}
\usepackage{amsmath}
\usepackage{amsfonts}
\usepackage{ulem}
\usepackage{graphicx,amsmath,amsfonts,amssymb,color,mathrsfs, amsmath,amsthm}
\usepackage{epsfig}
\newcommand{\E}{{\mathbb E}}

\newcommand{\chuhao}{\fontsize{19pt}{\baselineskip}\selectfont}

\newtheorem{thm}{Theorem}[section]
\newtheorem{prop}[thm]{Proposition}

\newtheorem{lem}[thm]{Lemma}
\newtheorem{defi}[thm]{Definition}
\newtheorem{remark}[thm]{Remark}
\newtheorem{example}[thm]{Example}
\newtheorem{pb}[thm]{Problem}

\newcommand{\be}{\begin{eqnarray*}}
\newcommand{\ee}{\end{eqnarray*}}
\newcommand{\beq}{\begin{equation}}
\newcommand{\eeq}{\end{equation}}

\numberwithin{equation}{section}

\title{\bf\color{black} \chuhao{The prdual and John-Nirenberg inequalities on generalized BMO martingale spaces}}
\author{Yong Jiao,\ \ Anming Yang,\ \ Lian Wu and Rui Yi
\\ { Central South University}
}

\date{\small 20 August, 2014}

\begin{document}
\maketitle


\makeatletter
\renewcommand{\@makefntext}[1]{#1}
\makeatother \footnotetext{\noindent
Yong Jiao is supported by NSFC (11471337), Hunan Provincial Natural Science Foundation(14JJ1004) 
and The International Postdoctoral Exchange Fellowship Program. 
\\ {\it Corresponding email:} jiaoyong@csu.edu.cn}


 \begin{abstract}
 In this paper we introduce the generalized BMO martingale spaces by stopping time
 sequences,  which enable us to characterize the dual spaces of martingale
 Hardy-Lorentz spaces $H_{p,q}^s$ for $0<p\leq1, 1<q<\infty$. Moreover,
 by duality we obtain a John-Nirenberg theorem for the generalized
 BMO martingale spaces when the stochastic basis is regular. We also extend the boundedness
 of fractional integrals to martingale Hardy-Lorentz spaces.
\end{abstract}


\section{Introduction}


\quad \quad  Basing mainly on the duality, John-Nirenberg inequality
and something else, the space BMO (Bounded Mean Oscillation; see
\cite{F1}, \cite{F2} and \cite{JN}) played a remarkable role in
classical analysis and probability. We refer to the monograghs \cite{BS}
and \cite{S2} for the function space version, respectively to the
monograghs \cite{B}, \cite{G1} and \cite{P} for the martingale version of
those theorems.

This paper deals with the John-Nirenberg inequalities and dualities
in the martingale theory. Before describing our main results, we
recall the classical John-Nirenberg inequalities in the martingale
theory. Let $(\Omega,\mathcal{F},\mathbb{P})$ be a complete
probability space and $\{\mathcal{F}_n\}_{n\geq0}$ be a
nondecreasing sequence of sub-$\sigma$-algebras of $\mathcal{F}$
such that $\mathcal{F}=\sigma(\bigcup\limits_n\mathcal{F}_n)$. We
also call $\{\mathcal{F}_n\}_{n\geq0}$ a stochastic basis
 (with convention $\mathcal{F}_{-1}=\mathcal{F}_0$).
The expectation operator and the conditional expectation operators
relative to $\mathcal{F}_n$ are denoted by $\mathbb{E}$ and
$\mathbb{E}_n$, respectively. The stochastic basis $\{\mathcal{F}_n\}_{n\geq0}$ is said to be
regular, if there exist an absolute constant $R>0$ such that
\begin{equation}
f_n\leq Rf_{n-1},
\end{equation}
holds for all nonnegative martingales $f=(f_n)_{n\geq0}$. 

A sequence $f = (f_n)_{n\geq0}$ of random variables such that $f_n$
is $\mathcal{F}_n$-measurable is said to be a martingale if $\E(|
f_n|)<\infty$ and $\E_n(f_{n+1}) = f_n$ for every $n\geq0.$ For the
sake of simplicity, we assume $f_0=0$. For $1\leq r<\infty$, the
Banach spaces $BMO_r$ are defined as follows:
$$BMO_r=\big\{ f= (f_n)_{n\geq0}\in L_r:\|f\|_{BMO_r}=\sup_n\|(\E_n|f-f_{n-1}|^r)^{\frac{1}{r}}\|_{\infty}<\infty \big\}.$$
Here $f$ in $|f-f_{n-1}|^r$ means $f_\infty$. The usual $BMO$
norm corresponds to $r=2$ above, i.e., $\|f\|_{BMO}=\|f\|_{BMO_2}$.
The John-Nirenberg theorem says that in the sense of equivalent
norms,
\begin{equation}
BMO_r=BMO, \ \ \ (1\leq r<\infty).
\end{equation}
A duality argument yields that (1.1) can be rewritten as follows
\begin{equation}
\|f\|_{BMO}\approx
\sup\limits_n\sup\limits_{A\in\mathcal{F}_n}\mathbb{P}(A)^{-\frac{1}
{r}}\Big(\int_A|f-f_{n-1}|^rd\mathbb{P}\Big)^{\frac{1}{r}}.
\end{equation}
Here and in the sequel, $ A\approx B$ means that there exist two
absolute constants $C_1$ and $C_2$ such that $C_1B\leq A\leq C_2B.$

The special contribution of this paper is to define the following
generalized BMO martingale space $BMO_{r,q}(\alpha)$ by stopping
time sequences.

\begin{defi}
For $1\leq r,q<\infty, \alpha\geq0$, the generalized BMO martingale
space is defined by
$$BMO_{r,q}(\alpha)=\Bigg\{f\in L_r: \|f\|_{BMO_{r,q}(\alpha)}=\sup\frac{\sum\limits_{k\in\mathbb{Z}}2^k \mathbb{P}(\nu_k<\infty)^{1-\frac{1}{r}}\|f-f^{\nu_k}\|_r}{\Big(\sum\limits_{k\in\mathbb{Z}}\big(2^k\mathbb{P}(\nu_k<\infty)^{1+\alpha}\big)^q\Big)^{\frac{1}{q}}}<\infty\Bigg\},$$
where the supremum is taken over all stopping time sequences
$\{\nu_k\}_{k\in\mathbb{Z}}$ such that
$\big\{2^k\mathbb{P}(\nu_k<\infty)^{1+\alpha}\big\}_{k\in\mathbb{Z}}\in
\ell_q$.
\end{defi}

Then the
generalized John-Nirenberg theorem, one of our main results, reads
as follows.
\begin{thm}
Suppose that the stochastic basis $\{\mathcal{F}_n\}_{n\geq0}$ is
regular and $1\leq q<\infty$. Then
\begin{equation}BMO_{r,q}(\alpha)=BMO_{2,q}(\alpha),
\end{equation} in the sense of equivalent norms for all $1\leq r<\infty$.
\end{thm}
We now explain the relation between (1.1) and (1.4). Let
$\mathcal{T}$ be the set of all stopping times relative to
$\{\mathcal{F}_n\}_{n\geq0}$.  On one hand, if the stopping time
sequence $\{\nu_k\}_{k\in\mathbb{Z}}$ reduces to a sequence whose
one element is a stopping time $\nu$ and the others are $\infty$,
then the generalized BMO space $BMO_{r,q}(\alpha)$ reduces to the
following Lipschitz space
$$BMO_r(\alpha)=\{f\in L_r: \|f\|_{BMO_r(\alpha)}=\sup\limits_{\nu\in\mathcal{T}}\mathbb{P}
(\nu<\infty)^{-\frac{1}{r}-\alpha}\|f-f^\nu\|_r<\infty\}.$$ On the
other hand, if the stochastic basis $\{\mathcal{F}_n\}_{n\geq0}$ is
regular, it is not very difficult to check that (1.2) can further be
reformulated as
\begin{equation}
\|f\|_{BMO}\approx \sup\limits_{\nu\in\mathcal{T}}\mathbb{P}
(\nu<\infty)^{-\frac{1}{r}}\|f-f^\nu\|_r\,\,.
\end{equation}
See also \cite{W2} for the facts above. Hence if $\alpha=0$, (1.4)
exactly implies (1.5). Consequently,  (1.1) can be deduced from (1.4)
when the stochastic basis $\{\mathcal{F}_n\}_{n\geq0}$ is regular.

We now turn to the second aim of this paper. The generalized $BMO$
martingale space $BMO_{r,q}(\alpha)$ defined in this paper enable us
to characterize the dualities of martingale Hardy-Lorentz spaces for
$0<p\leq1, 1<q<\infty$. It is well known that the dual spaces of
Lebesgue spaces $L_p$
 or Lorentz spaces $L_{p,q}$ are trivial when $0<p<1$ (see for instance \cite{C} or \cite{G2}), namely,
$$\ \big(L_p\big)^*=\big(L_{p,q}\big)^*=\{0\}, \ \ \ \ \ (0<p<1,\,0<q\leq\infty).$$
However, the dual spaces of martingale Hardy spaces are very
different from those of Lebesgue spaces $L_p$ and Lorentz spaces
$L_{p,q}$. This can be illustrated by the fact that the dual spaces
of $L_p$ and $L_{p,q}$ $(0<p<1)$ are trivial while
$$\big(H_p^s\big)^*=BMO_2(\alpha), \quad (0<p<1,\ \alpha=\frac{1}{p}-1),$$
where $H_p^s$ denotes the martingale Hardy space associated with the
conditional quadratic variation, that is,
$$H_{p}^s=\Big\{f=(f_n)_{n\geq0}: \|f\|_{H_{p}^s}=\Big\|
\Big(\sum\limits_{i=1}^\infty
\mathbb{E}_{i-1}|d_if|^2\Big)^{\frac{1}{2}}\Big\|_{p}<\infty\Big\}.$$
We refer to \cite{G1}, \cite{L} and \cite{W2} for the fact above. At
the same time, Weisz \cite{W2} also proved the following duality
result for martingale Hardy-Lorentz spaces,
$$\big(H_{p,q}^s\big)^*=H_{p',q'}^s, \quad (1<p<\infty,\ 1\leq q<\infty),$$
where $p'$ and $q'$ denote the conjugate numbers of $p$ and $q$
respectively; see Section 2 for definition of $H_{p,q}^s$. But the
question how to characterize the dual spaces of martingale
Hardy-Lorentz spaces for $0<p\leq1$, $0<q<\infty$ is still open.  We
prove that the dual space of martingale Hardy-Lorentz space is the
same as the one of martingale Hardy spaces when $0<p,q\leq1$, while
it needs the new notion $BMO_{r,q}(\alpha)$ when $0<p\leq1$,
$1<q<\infty$. In Section 4 we shall show
\begin{thm} Let $0<p\leq1,$ $\alpha=\frac{1}{p}-1.$ Then
$$\big(H_{p,q}^s\big)^*=BMO_2(\alpha),\quad 0<q\leq1;$$ and
$$\big(H_{p,q}^s\big)^*=BMO_{2,q}(\alpha), \quad 1<q<\infty.$$
\end{thm}

This paper will be divided into five further sections. In the next
section, some notations and basic knowledge will be introduced. In
Section 3, the atomic decompositions of martingale Hardy-Lorentz
spaces are formulated. In Section 4, using atomic decompositions in Section 3, we prove some dual theorems of
martingale Hardy-Lorentz spaces. By duality, the new John-Nirenberg
theorem for the generalized BMO martingale space is proved in Section 5.
In the final Section, the boundedness of fractional integrals on
martingale Hardy-Lorentz spaces are investigated.

In this paper, the set of integers and the set of nonnegative
integers are always denoted by $\mathbb{Z}$ and $\mathbb{N}$,
respectively. We use $C$ to denote the absolute constant which may
vary from line to line.


\section{Notations and preliminaries}


We first introduce the distribution function and the decreasing
rearrangement. Let $f$ be a measurable function defined on the
probability space $(\Omega,\mathcal{F},\mathbb{P})$. We define the
distribution function of $f$ by
$$\lambda_s(f)=\mathbb{P}\big(\{\omega\in\Omega: |f(\omega)|>s\}\big), \ \ (s\geq0).$$
And denote by $\mu_t(f)$ the decreasing rearrangement of $f$,
defined by
$$\mu_t(f)=\inf\{s\geq0: \lambda_s(f)\leq t\}, \ \ (t\geq0),$$ with the convention that $\inf \emptyset=\infty$.

We list some properties of distribution functions and decreasing
rearrangements in the following proposition. The properties will be
used in the proof of theorems in the later sections.

\begin{prop}
Let $f$ and $g$ be two measurable functions on
$(\Omega,\mathcal{F},\mathbb{P})$, then we have

(1) if $|f|\leq |g|$ $\mathbb{P}$-$a.e.$ then
$\lambda_s(f)\leq\lambda_s(g)$ for all $s\geq0$;

(2) $\lambda_{s_1+s_2}(f+g)\leq\lambda_{s_1}(f)+\lambda_{s_2}(g)$
for all $s_1, s_2\geq0$;

(3) $\mu_t(af)=|a|\mu_t(f)$ for all $a\in \mathbb{C}$ and $t\geq0$;

(4) if $|f|\leq |g|$ $\mathbb{P}$-$a.e.$ then $\mu_t(f)\leq\mu_t(g)$
for all $t\geq0$;

(5) $\mu_{t_1+t_2}(f+g)\leq\mu_{t_1}(f)+\mu_{t_2}(g)$ for all $t_1,
t_2\geq0$.
\end{prop}

The Lorentz space $L_{p,q}(\Omega,\mathcal{F},\mathbb{P}),\,
0<p<\infty, 0<q\leq\infty$, consists of those measurable functions
$f$ with finite norm or quasinorm $\|f\|_{p,q}$ given by
$$\|f\|_{p,q}=\Bigg(\frac{q}{p}\int_0^\infty\Big(t^{\frac{1}{p}}\mu_t(f)\Big)^q\frac{dt}{t}\Bigg)^{\frac{1}{q}}, \ \ (0<q<\infty),$$
$$\|f\|_{p,\infty}=\sup\limits_{t>0}t^{\frac{1}{p}}\mu_t(f), \ \ (q=\infty).$$

It will be convenient for us to use an equivalent definition of
$\|f\|_{p,q}$, namely
$$\|f\|_{p,q}=\Bigg(q\int_0^\infty\bigg(t\mathbb{P}(|f(x)|>t)^{\frac{1}{p}}\bigg)^q\frac{dt}{t}\Bigg)^{\frac{1}{q}}, \ \ (0<q<\infty),$$
$$\|f\|_{p,\infty}=\sup\limits_{t>0}t\mathbb{P}(|f(x)|>t)^{\frac{1}{p}}, \ \ (q=\infty).$$

We recall that Lorentz spaces $L_{p,q}$ increase as the second
exponent $q$ increases, and decrease as the first exponent $p$
increases (the second exponent $q$ is not involved). Namely,
$L_{p,q_1}\subset L_{p,q_2}$ for $0<p<\infty$ and $0<q_1\leq q_2\leq
\infty$, $L_{p_1,q_1}\subset L_{p_2,q_2}$ for $0<p_2<p_1<\infty$ and
$0<q_1,q_2\leq\infty$. It is also well known that if $1<p<\infty$
and $1\leq q\leq\infty$, or $p=q=1$, then $\|\cdot\|_{p,q}$ is
equivalent to a norm. However, for the other values of $p$ and $q$,
$\|\cdot\|_{p,q}$ is only a quasi-norm. In particular, if $0<q\leq1$
and $q\leq p<\infty$, then $\|\cdot\|_{p,q}$ is equivalent to a
$q$-norm. H\"{o}lder's inequality for Lorentz spaces is the
following
$$\|fg\|_{p,q}\leq C\|f\|_{p_1,q_1}\|g\|_{p_2,q_2},$$
where $0<p,p_1,p_2<\infty$ and $0<q,q_1,q_2\leq\infty$ such that
$\frac{1}{p}=\frac{1}{p_1}+\frac{1}{p_2}$ and
$\frac{1}{q}=\frac{1}{q_1}+\frac{1}{q_2}$.

We now introduce martingale Hardy-Lorentz spaces. Denote by
$\mathcal{M}$ the set of all martingales $f=(f_n)_{n\geq0}$ relative
to $\{\mathcal{F}_n\}_{n\geq0}$ such that $f_0=0$. For
$f\in\mathcal{M}$, denote its martingale difference by
$d_nf=f_n-f_{n-1}$ ($n\geq0$, with convention $f_{-1}=0)$. Then the
maximal function, the quadratic variation and the conditional
quadratic variation of a martingale $f$ are respectively defined by
$$f_n^*=\sup\limits_{0\leq i\leq n}|f_i|, \ \ \
f^*=\sup\limits_{n\geq0}|f_n|,$$
$$S_n(f)=\Big(\sum\limits_{i=1}^n|d_if|^2\Big)^{\frac{1}{2}},\ \ \ S(f)=\Big(\sum\limits_{i=1}^\infty|d_if|^2\Big)^{\frac{1}{2}},$$
$$s_n(f)=\Big(\sum\limits_{i=1}^n\mathbb{E}_{i-1}|d_if|^2\Big)^{\frac{1}{2}},\ \ \ s(f)=\Big(\sum\limits_{i=1}^\infty \mathbb{E}_{i-1}|d_if|^2\Big)^{\frac{1}{2}}.$$
Let $\Lambda$ be the collection of all sequences
$(\lambda_n)_{n\geq0}$ of nondecreasing, nonnegative and adapted
functions, set
$\lambda_\infty=\lim\limits_{n\rightarrow\infty}\lambda_n$. For
$f\in\mathcal{M}, 0<p<\infty, 0<q\leq\infty$, let
$$\Lambda[Q_{p,q}](f)=\{(\lambda_n)_{n\geq0}\in\Lambda: S_n(f)\leq\lambda_{n-1}(n\geq1), \lambda_\infty\in L_{p,q}\},$$
$$\Lambda[D_{p,q}](f)=\{(\lambda_n)_{n\geq0}\in\Lambda: |f_n|\leq\lambda_{n-1}(n\geq1), \lambda_\infty\in L_{p,q}\}.$$
We define martingale Hardy-Lorentz spaces as follows. For
$0<p<\infty, 0<q\leq\infty$,
$$H_{p,q}^*=\{f\in\mathcal{M}: \|f\|_{H_{p,q}^*}=\|f^*\|_{p,q}<\infty\},$$
$$H_{p,q}^S=\{f\in\mathcal{M}: \|f\|_{H_{p,q}^S}=\|S(f)\|_{p,q}<\infty\},$$
$$H_{p,q}^s=\{f\in\mathcal{M}: \|f\|_{H_{p,q}^s}=\|s(f)\|_{p,q}<\infty\},$$
$$Q_{p,q}=\{f\in\mathcal{M}: \|f\|_{Q_{p,q}}=\inf\limits_{(\lambda_n)_{n\geq0}\in\Lambda[Q_{p,q}](f)}\|\lambda_\infty\|_{p,q}<\infty\},$$
$$D_{p,q}=\{f\in\mathcal{M}: \|f\|_{D_{p,q}}=\inf\limits_{(\lambda_n)_{n\geq0}\in\Lambda[D_{p,q}](f)}\|\lambda_\infty\|_{p,q}<\infty\}.$$
If taking $p=q$ in the definitions above, we get the usual
martingale Hardy spaces. In order to describe the duality theorems,
we need to introduce the Lipschitz space $BMO_r(\alpha)$. For $1\leq
r<\infty$, $\alpha\geq 0$, the Lipschitz space are defined as
follows
$$BMO_r(\alpha)=\{f\in L_r: \|f\|_{BMO_r(\alpha)}<\infty\},$$
where
$$\|f\|_{BMO_r(\alpha)}=\sup\limits_{n\in\mathbb{N}}\sup\limits_{A\in\mathcal{F}_n}\mathbb{P}(A)^{-\frac{1}
{r}-\alpha}\Big(\int_A|f-\mathbb{E}_nf|^rd\mathbb{P}\Big)^{\frac{1}{r}}.$$
\noindent Let $\mathcal{T}$ be the set of all stopping times
relative to $\{\mathcal{F}_n\}_{n\geq0}$. For a martingale
$f=(f_n)_{n\geq0}\in\mathcal{M}$ and a stopping time
$\nu\in\mathcal{T}$, we denote the stopped martingale by
$f^\nu=(f_n^\nu)_{n\geq0}=(f_{n\wedge\nu})_{n\geq0}$. Then it is
easy to show that
$$\|f\|_{BMO_r(\alpha)}=\sup\limits_{\nu\in\mathcal{T}}\mathbb{P}(\nu<\infty)^{-\frac{1}{r}-\alpha}\|f-f^\nu\|_r, \ \ \ (1\leq r<\infty, \alpha\geq0).$$

The main new notion of the present paper is the generalized BMO
martingale space $BMO_{r,q}(\alpha)$ $(1\leq r,q<\infty, \
\alpha\geq0)$, see Section 1 for the definition. In Definition 1.1,
if the stopping time sequence $\{\nu_k\}_{k\in\mathbb{Z}}$ reduces
to a sequence whose one element is a stopping time $\nu$ and the
others are $\infty$, then the generalized BMO martingale space
$BMO_{r,q}(\alpha)$ reduces to the Lipschitz martingale space $BMO_r(\alpha)$.
Obviously, $BMO_{r,q}(\alpha)$ is a subspace of $BMO_r(\alpha)$ and
$\|f\|_{BMO_r(\alpha)}\leq\|f\|_{BMO_{r,q}(\alpha)}$.

We will present the atomic decomposition theorems for martingale
Hardy-Lorentz spaces in the next section. Now let us introduce the
notion of atoms; see for example \cite{W2}.

\begin{defi}
A measurable function $a$ is called a $(p,\infty)$-atom of the first
category (or of the second category, or of the third category) if
there exists a stopping time $\nu\in\mathcal{T}$ ($\nu$ is called
the stopping time associated with a) such that

(i) $a_n=\mathbb{E}_n(a)=0$, (if $\nu\geq n$), (ii)
$\|s(a)\|_\infty\leq \mathbb{P}(\nu<\infty)^{-\frac{1}{p}}$ (or
(ii$'$) $\|S(a)\|_\infty\leq \mathbb{P}(\nu<\infty)^{-\frac{1}{p}}$,
or (ii$''$) $\|a^*\|_\infty\leq
\mathbb{P}(\nu<\infty)^{-\frac{1}{p}}$, respectively).
\end{defi}\noindent These three category atoms are briefly called $(1,p,\infty)$-atom,
$(2,p,\infty)$-atom and $(3,p,\infty)$-atom, respectively.

We conclude this section by two lemmas which are very useful to
verify that a function is in Lorentz spaces $L_{p,q}$, which are respectively from
 Lemma 1.1 and Lemma 1.2 in \cite{AT}.

\begin{lem}
Let $0<p<\infty, 0<q\leq\infty$, assume that the nonnegative
sequence $\{\mu_k\}$ satisfies $\{2^k\mu_k\}\in l^q$. Further
suppose that the nonnegative function $\varphi$ verifies the
following property: there exists $0<\varepsilon<1$ such that, given
an arbitrary integer $k_0$, we have
$\varphi\leq\psi_{k_0}+\eta_{k_0}$, where $\psi_{k_0}$ is
essentially bounded and satisfies $\|\psi_{k_0}\|_\infty\leq
C2^{k_0}$, and
$$2^{k_0\varepsilon p}\mathbb{P}(\eta_{k_0}>2^{k_0})\leq C\sum\limits_{k=k_0}^\infty(2^{k\varepsilon}\mu_k)^p.$$
Then $\varphi\in L_{p,q}$ and $\|\varphi\|_{p,q}\leq
C\|\{2^k\mu_k\}\|_{l_q}$.
\end{lem}

\begin{lem}
Let $0<p<\infty$, and let the nonnegative sequence $\{\mu_k\}$ be
such that $\{2^k\mu_k\}\in l^q$, $0<q\leq\infty$. Further, suppose
that the nonnegative function $\varphi$ satisfies the following
property: there exists $0<\varepsilon<1$ such that, given an
arbitrary integer $k_0$, we have $\varphi\leq\psi_{k_0}+\eta_{k_0}$,
where $\psi_{k_0}$ and $\eta_{k_0}$ satisfy
$$2^{k_0p}\mathbb{P}(\psi_{k_0}>2^{k_0})^\varepsilon\leq C\sum\limits_{k=-\infty}^{k_0}(2^k\mu_k^\varepsilon)^p, \ \ \ 0<\varepsilon<\min(1,\frac{q}{p}),$$
$$2^{k_0\varepsilon p}\mathbb{P}(\eta_{k_0}>2^{k_0})\leq C\sum\limits_{k=k_0}^\infty(2^{k\varepsilon}\mu_k)^p.$$
Then $\varphi\in L_{p,q}$ and $\|\varphi\|_{p,q}\leq
C\|\{2^k\mu_k\}\|_{l_q}$.
\end{lem}


\section{Atomic decompositions}


\quad \quad The method of atomic decompositions plays an important
role in martingale theory; see for instance \cite{HR}, \cite{LH},
\cite{ML}, \cite{MNS}, \cite{W1} and \cite{W2}. In particular, Jiao, Peng and Liu \cite{JPL} proved the atomic decompositions
of martingale Hardy-Lorentz spaces in 2009. Since $\|\cdot\|_{p,q}$ is
equivalent to a $q$-norm just when $0<q\leq1$ and $q\leq p<\infty$,
there is a restrictive condition for the converse part of Theorem
2.1 in \cite{JPL}. We improve Theorem 2.1 in \cite{JPL} by using the
technical Lemma 2.3 and shows that the converse part of Theorem 2.1
in \cite{JPL} is true for all $0<p<\infty, 0<q\leq\infty$.

\begin{thm}
If $f=(f_n)_{n\geq0}\in H_{p,q}^s $ $(0<p<\infty, 0<q\leq\infty)$,
then there exists a sequence $(a^k)_{k\in\mathbb{Z}}$ of
$(1,p,\infty)$-atoms and a sequence $(\mu_k)_{k\in\mathbb{Z}}\in
l_q$ of real numbers satisfying
$\mu_k=A\cdot2^k\mathbb{P}(\nu_k<\infty)^{\frac{1}{p}}$ (where $A$
is a positive constant and $\nu_k$ is the stopping time associated
with $a^k$) such that \begin{eqnarray}f_n=\sum\limits_{k\in
\mathbb{Z}}\mu_ka_n^k,\quad a.e., \ \ \ n\in
\mathbb{N},\end{eqnarray} and
$$\|(\mu_k)_{k\in\mathbb{Z}}\|_{l_q}\leq C\|f\|_{H_{p,q}^s}.$$
Conversely, if the martingale $f$ has the above decomposition, then
$f\in H_{p,q}^s$ and
$$\|f\|_{H_{p,q}^s}\approx\inf\|(\mu_k)_{k\in\mathbb{Z}}\|_{l_q},$$
where the infimum is taken over all the above decompositions.
\end{thm}

{Proof.} Assume that $f\in H_{p,q}^s$ $(0<p<\infty,0<q\leq\infty)$.
For each $k\in\mathbb{Z}$, the stopping time is defined as follows
$$\nu_k=\inf\{n\in\mathbb{N}: s_{n+1}(f)>2^k\}, \ \ \ (\inf\emptyset=\infty).$$
Obviously, the sequence of these stopping times is non-decreasing.
Similarly to the proof of Theorem 2.2 in \cite{W2} ( or see the
proof of Theorem 2.1 in \cite{JPL}), we have
\begin{eqnarray*}
\sum\limits_{k\in\mathbb{Z}}(f_n^{\nu_{k+1}}-f_n^{\nu_k})=f_n.
\end{eqnarray*}
Let
$$\mu_k=3\cdot2^k\mathbb{P}(\nu_k<\infty)^{\frac{1}{p}},\ \ \ a_n^k=\frac{f_n^{\nu_{k+1}}-f_n^{\nu_k}}{\mu_k}.$$
If $\mu_k=0$, we assume that $a_n^k=0$. Then for any fixed
$k\in\mathbb{Z}$, $a^k=(a_n^k)_{n\geq0}$ is a martingale.
Considering the stopped martingale
$f^{\nu_k}=(f_n^{\nu_k})_{n\geq0}=(f_{n\wedge\nu_k})_{n\geq0}$, we
have $s(f^{\nu_k})=s_{\nu_k}(f)\leq2^k$,
$s(f^{\nu_{k+1}})\leq2^{k+1}$. Then
$$s(a^k)\leq\frac{s(f^{\nu_{k+1}})+s(f^{\nu_k})}{\mu_k}\leq\mathbb{P}(\nu_k<\infty)^{-\frac{1}{p}},$$
which implies that $(a_n^k)_{n\geq0}$ is a $L_2$-bounded martingale.
So $(a_n^k)_{n\geq0}$ converges in $L_2$. Denote the limit still by
$a^k$, then $\mathbb{E}_na^k=a_n^k$. If $\nu_k\geq n$, then
$a_n^k=0$, and
$\|s(a^k)\|_\infty\leq\mathbb{P}(\nu_k<\infty)^{-\frac{1}{p}}$. Thus
we conclude that $a^k$ is really a $(1,p,\infty)$-atom. Since
$\{\nu_k<\infty\}=\{s(f)>2^k\}$, we get for $0<q<\infty$
\begin{eqnarray*}
\Big(\sum\limits_{k\in\mathbb{Z}}|\mu_k|^q\Big)^{\frac{1}{q}}&=&3\Big(\sum\limits_{k\in\mathbb{Z}}2^{kq}\mathbb{P}(\nu_k<\infty)^{\frac{q}{p}}\Big)^{\frac{1}{q}}
=3\Big(\sum\limits_{k\in\mathbb{Z}}2^{kq}\mathbb{P}(s(f)>2^k)^{\frac{q}{p}}\Big)^{\frac{1}{q}}
\\&\leq&C\Big(\sum\limits_{k\in\mathbb{Z}}\int_{2^{k-1}}^{2^k}y^{q-1}dy \mathbb{P}(s(f)>2^k)^{\frac{q}{p}}\Big)^{\frac{1}{q}}
\\&\leq&C\Big(\sum\limits_{k\in\mathbb{Z}}\int_{2^{k-1}}^{2^k}y^{q-1}\mathbb{P}(s(f)>y)^{\frac{q}{p}}dy\Big)^{\frac{1}{q}}
\\&\leq&C\Big(\int_0^\infty y^{q-1}\mathbb{P}(s(f)>y)^{\frac{q}{p}}dy\Big)^{\frac{1}{q}}
\\&\leq&C\|s(f)\|_{p,q}=C\|f\|_{H_{p,q}^s}.
\end{eqnarray*}
For $q=\infty,$
\begin{eqnarray*}
\|(\mu_k)_{k\in\mathbb{Z}}\|_\infty&=&\sup\limits_{k\in\mathbb{Z}}|\mu_k|=3\cdot\sup\limits_{k\in\mathbb{Z}}2^k
\mathbb{P}(\nu_k<\infty)^{\frac{1}{p}}
\\&=&3\cdot\sup\limits_{k\in\mathbb{Z}}2^k \mathbb{P}(s(f)>2^k)^{\frac{1}{p}}
\\&\leq&C\|s(f)\|_{p,\infty}=C\|f\|_{H_{p,\infty}^s}.
\end{eqnarray*}
Consequently, $\|(\mu_k)_{k\in\mathbb{Z}}\|_{l_q}\leq
C\|f\|_{H_{p,q}^s}$.

Conversely, if the martingale $f$ has the above decomposition, then
for an arbitrary integer $k_0$, let
$$f_n=\sum\limits_{k\in\mathbb{Z}}\mu_ka_n^k=g_n+h_n, \ \ \ (n\in\mathbb{N}),$$
where $g_n=\sum\limits_{k=-\infty}^{k_0-1}\mu_ka_n^k$ and
$h_n=\sum\limits_{k=k_0}^\infty\mu_ka_n^k$. By the sublinearity of
the operator $s$, we have $s(f)\leq s(g)+s(h)$. Then
\begin{eqnarray*}
\|s(g)\|_\infty&\leq&\|\sum\limits_{k=-\infty}^{k_0-1}|\mu_k|s(a^k)\|_\infty\leq\sum\limits_{k=-\infty}^{k_0-1}|\mu_k|\|s(a^k)\|_\infty
\\&\leq&\sum\limits_{k=-\infty}^{k_0-1}|\mu_k|\mathbb{P}(\nu_k<\infty)^{-\frac{1}{p}}
\\&\leq&\sum\limits_{k=-\infty}^{k_0-1}A\cdot2^k=A\cdot2^{k_0}.
\end{eqnarray*}
Since $s(a^k)=0$ on the set $\{\nu_k=\infty\}$, we have
$\{s(a^k)>0\}\subset\{\nu_k<\infty\}$. Then it follows from
$s(h)\leq\sum\limits_{k=k_0}^\infty|\mu_k|s(a^k)$ that
$$\{s(h)>0\}\subset\bigcup\limits_{k=k_0}^\infty\{s(a^k)>0\}\subset\bigcup\limits_{k=k_0}^\infty\{\nu_k<\infty\}.$$
The for each $0<\varepsilon<1$, we obtain
\begin{eqnarray*}
2^{k_0\varepsilon p}\mathbb{P}(s(h)>2^{k_0})&\leq&2^{k_0\varepsilon
p}\mathbb{P}(s(h)>0)\leq2^{k_0\varepsilon
p}\sum\limits_{k=k_0}^\infty \mathbb{P}(\nu_k<\infty)
\\&=&2^{k_0\varepsilon p}\sum\limits_{k=k_0}^\infty2^{k\varepsilon p}\mathbb{P}(\nu_k<\infty)2^{-k\varepsilon p}
\\&\leq&\sum\limits_{k=k_0}^\infty2^{k\varepsilon p}\mathbb{P}(\nu_k<\infty)=\sum\limits_{k=k_0}^\infty\Big(2^{k\varepsilon}\mathbb{P}(\nu_k<\infty)^{\frac{1}{p}}\Big)^p.
\end{eqnarray*}
By Lemma 2.3, we have $s(f)\in L_{p,q}$ and $\|s(f)\|_{p,q}\leq
C\|\{2^k\mathbb{P}(\nu_k<\infty)^{\frac{1}{p}}\}_{k\in\mathbb{Z}}\|_{l_q}\leq
C\|(\mu_k)_{k\in\mathbb{Z}}\|_{l_q}$. Then $f\in H_{p,q}^s$ and
$\|f\|_{H_{p,q}^s}\leq C\|(\mu_k)_{k\in\mathbb{Z}}\|_{l_q}$. Thus
$$\|f\|_{H_{p,q}^s}\approx\inf\|(\mu_k)_{k\in\mathbb{Z}}\|_{l_q},$$
where the infimum is taken over all the above decompositions. The
proof of the theorem is complete.

\begin{remark}
 If $q\neq\infty$, then (3.1) holds in $H_{p,q}^s$. Namely,
 the sum $\sum_{k=m}^n\mu_k a^k$ converges to $f$ in $H_{p,q}^s$
as $m\rightarrow -\infty$, $n\rightarrow \infty$.
 Indeed,
 $$\sum_{k=m}^n \mu_ka^k = \sum_{k=m}^n (f^{\nu_{k+1}}-f^{\nu_k})=f^{\nu_{n+1}}-f^{\nu_m}.$$
 By the sublinearity of $s(f)$ we have
 \be
 \|f-\sum_{k=m}^n\mu_k a^k\|_{H_{p,q}^s}&=&\|s(f-f^{\nu_{n+1}}+f^{\nu_m})\|_{{p,q}}\leq
 \|s(f-f^{\nu_{n+1}})+s(f^{\nu_m})\|_{{p,q}}\\
 &\leq & C
 \bigg(\|s(f-f^{\nu_{n+1}})\|_{{p,q}}+\|s(f^{\nu_m})\|_{{p,q}}\bigg).
 \ee
 Since $s(f-f^{\nu_{n+1}})^2=s(f)^2- s(f^{\nu_{n+1}})^2$, then
 $s(f-f^{\nu_{n+1}})\leq s(f)$, $s(f^{\nu_m})\leq s(f)$
 and $s(f-f^{\nu_{n+1}}),\; s(f^{\nu_m})\rightarrow 0 $ a.e. as $m\rightarrow -\infty$, $n\rightarrow \infty$.
 Thus by the Lebesgue convergence theorem, we have
 $$\|s(f-f^{\nu_{n+1}})\|_{{p,q}},\;\|s(f^{\nu_m})\|_{{p,q}} \rightarrow 0\;\;\;\; as\; m\rightarrow -\infty,\;n\rightarrow \infty,$$
 which means $\|f-\sum_{k=m}^n\mu_k a^k\|_{H_{p,q}^s}\rightarrow 0$ as $m\rightarrow -\infty$, $n\rightarrow
 \infty$. Further, for each $k\in\mathbb {Z}$, $a^k=(a_n^k)_{n\geq 0}$ is $L_2$
 bounded, hence $H_2^s=L_2$ is dense in $H_{p,q}^s.$
 \end{remark}

\begin{thm}
In Theorem 3.1, if we  replace $H_{p,q}^s, (1,p,\infty)$-atoms by
$Q_{p,q}, (2,p,\infty)$-atoms (or $D_{p,q}, (3,p,\infty)$-atoms)
respectively, then the conclusions still hold.
\end{thm}
{ Proof.} The proof is similar to the one of Theorem 3.1, so we give
it in sketch, only. If $f=(f_n)_{n\geq0}\in Q_{p,q}\ (or \
D_{p,q})$. The stopping times $\nu_k$ are defined in these cases by
$$\nu_k=\inf\{n\in\mathbb{N}: \lambda_n>2^k\}, \ \ \ (\inf\emptyset=\infty),$$
where $(\lambda_n)_{n\geq0}$ is the sequence in the definition of
$Q_{p,q}\ (or \ D_{p,q})$. Let $a^k$ and $\mu_k$ $(k\in\mathbb{Z})$
be defined as in the proof of Theorem 3.1. Then the conclusions
$f_n=\sum\limits_{k\in\mathbb{Z}}\mu_ka_n^k$ $(n\in\mathbb{N})$ and
$\|(\mu_k)_{k\in\mathbb{Z}}\|_{l_q}\leq C\|f\|_{Q_{p,q}}\ (or \
\|(\mu_k)_{k\in\mathbb{Z}}\|_{l_q}\leq C\|f\|_{D_{p,q}})$ still
hold.

To prove the converse part, let
$$\lambda_n=\sum\limits_{k\in\mathbb{Z}}\mu_k\chi_{\{\nu_k\leq n\}}\|S(a^k)\|_\infty\ (or\  \lambda_n=\sum\limits_{k\in\mathbb{Z}}\mu_k\chi_{\{\nu_k\leq n\}}\|(a^k)^*\|_\infty).$$
Then $(\lambda_n)_{n\geq0}$ is a nondecreasing, nonnegative and
adapted sequence with $S_{n+1}(f)\leq\lambda_n\ (or \
|f_{n+1}|\leq\lambda_n)$.

For any given integer $k_0$, let
$$\lambda_\infty=\lambda_\infty^{(1)}+\lambda_\infty^{(2)},$$ where
$$\lambda_\infty^{(1)}=\sum\limits_{k=\infty}^{k_0-1}\chi_{\{\nu_k<\infty\}}\|S(a^k)\|_\infty\ (or\ \lambda_\infty^{(1)}=\sum\limits_{k=\infty}^{k_0-1}\chi_{\{\nu_k<\infty\}}\|(a^k)^*\|_\infty),$$
and
$$\lambda_\infty^{(2)}=\sum\limits_{k=k_0}^\infty\chi_{\{\nu_k<\infty\}}\|S(a^k)\|_\infty\ (or\  \lambda_\infty^{(2)}=\sum\limits_{k=k_0}^\infty\chi_{\{\nu_k<\infty\}}\|(a^k)^*\|_\infty).$$
Replacing $s(g)$ and $s(h)$ in the proof of Theorem 3.1 by
$\lambda_\infty^{(1)}$ and $\lambda_\infty^{(2)}$. Using Lemma 2.3,
we can obtain $f\in Q_{p,q}\ (or \ f\in D_{p,q})$ and
$\|f\|_{Q_{p,q}}\approx\inf\|(\mu_k)_{k\in\mathbb{Z}}\|_{l_q}\ (or \
\|f\|_{D_{p,q}}\approx\inf\|(\mu_k)_{k\in\mathbb{Z}}\|_{l_q})$,
where the infimum is taken over all the above decompositions. The
proof is complete.


\section{Duality results}

\quad \quad In this section, we prove the predual of the generalized
BMO martingale spaces.
\begin{thm}
The dual space of $H_{p,q}^s$ is $BMO_2(\alpha)$, $(0<p,q\leq 1,
\alpha=\frac{1}{p}-1)$.
\end{thm}

{ Proof.} Since $0<p,q\leq1$, we note that by
$$\|f\|_{H_{p,q}^s}=\|s(f)\|_{p,q}\leq\|s(f)\|_{2,2}=\|f\|_2,$$
the space $L_2$ is a subspace of $H_{p,q}^s$. By the Remark 3.2, we
know that $L_2$ is dense in $H_{p,q}^s$. For any
$g\in BMO_2{(\alpha)}\subset L_2$, we show that
$$\varphi_g(f)=\mathbb{E}(fg), \ \ \ \forall f\in L_2,$$ is a continuous
linear functional on $L_2$. It follows from Theorem 3.1 that
$f=\sum\limits_{k\in \mathbb{Z}}\mu_k a^k$. Hence
$$\varphi_g(f)=\mathbb{E}(fg)=\sum\limits_{k\in\mathbb{Z}}\mu_k\mathbb{E}(a^kg).$$
By the definition of the atom $a^k$, we have
$\mathbb{E}(a^kg)=\mathbb{E}(a^k(g-g^{\nu_k}))$. Using H\"{o}lder's
inequality, we obtain
\begin{eqnarray*}
|\varphi_g(f)|&=&|\sum\limits_{k\in\mathbb{Z}}\mu_k\mathbb{E}(a^k(g-g^{\nu_k}))|\leq\sum\limits_{k\in\mathbb{Z}}|\mu_k|\mathbb{E}\big(|a^k(g-g^{\nu_k})|\big)
\\&\leq&\sum\limits_{k\in\mathbb{Z}}|\mu_k|\|a^k\|_2\|g-g^{\nu_k}\|_2=\sum\limits_{k\in\mathbb{Z}}|\mu_k|\|s(a^k)\|_2\|g-g^{\nu_k}\|_2
\\&=&\sum\limits_{k\in\mathbb{Z}}|\mu_k|\|s(a^k)\chi_{\{\nu_k<\infty\}}\|_2\|g-g^{\nu_k}\|_2
\leq\sum\limits_{k\in\mathbb{Z}}|\mu_k|\|s(a^k)\|_\infty\|\chi_{\{\nu_k<\infty\}}\|_2\|g-g^{\nu_k}\|_2
\\&\leq&\sum\limits_{k\in\mathbb{Z}}|\mu_k|\mathbb{P}(\nu_k<\infty)^{\frac{1}{2}-\frac{1}{p}}\|g-g^{\nu_k}\|_2
\leq\sum\limits_{k\in\mathbb{Z}}|\mu_k|\|g\|_{BMO_2(\alpha)}.
\end{eqnarray*}
Since $0<q\leq1$, we have
$|\varphi_g(f)|\leq\Big(\sum\limits_{k\in\mathbb{Z}}|\mu_k|^q\Big)^{\frac{1}{q}}\|g\|_{BMO_2(\alpha)}$,
and by Theorem 3.1, we obtain
$$|\varphi_g(f)|\leq C\|f\|_{H_{p,q}^s}\|g\|_{BMO_2(\alpha)}.$$
By density of $L_2$ in $H_{p,q}^s$, $\varphi_g$ can be uniquely
extended to a continuous functional on $H_{p,q}^s$.

Conversely, for any $\varphi\in(H_{p,q}^s)^*$, we show that there
exists $g\in BMO_2(\alpha)$ such that $\varphi=\varphi_g$ and
$\|g\|_{BMO_2(\alpha)}\leq\|\varphi\|$.

Since $L_2$ can be continuously embedded in $H_{p,q}^s$, then there
exists $g\in L_2$ such that $\varphi(f)=\mathbb{E}(fg)$, $(f\in
L_2)$.

Let $\nu$ be an arbitrary stopping time and
$$h=\frac{g-g^\nu}{\|g-g^\nu\|_2\mathbb{P}(\nu<\infty)^{\frac{1}{p}-\frac{1}{2}}}.$$
Then $s(h)=0$ on $\{\nu=\infty\}$, namely,
$s(h)=s(h)\chi_{\{\nu<\infty\}}$.

Since $0<p,q\leq1$, then there exists $p_1,q_1>0$ such that
$\frac{1}{p}=\frac{1}{2}+\frac{1}{p_1},
\frac{1}{q}=\frac{1}{2}+\frac{1}{q_1}$.

By H\"{o}lder's inequality we have
\begin{eqnarray*}
\|h\|_{H_{p,q}^s}&=&\frac{\|g-g^\nu\|_{H_{p,q}^s}}{\|g-g^\nu\|_2\mathbb{P}(\nu<\infty)^{\frac{1}{p}-\frac{1}{2}}}
=\frac{\|s(g-g^\nu)\chi_{\{\nu<\infty\}}\|_{p,q}}{\|g-g^\nu\|_2\mathbb{P}(\nu<\infty)^{\frac{1}{p}-\frac{1}{2}}}
\\&\leq&\frac{C}{\|g-g^\nu\|_2\mathbb{P}(\nu<\infty)^{\frac{1}{p}-\frac{1}{2}}}\|s(g-g^\nu)\|_{2,2}\|\chi_{\{\nu<\infty\}}\|_{p_1,q_1}
\\&=&\frac{C\|g-g^\nu\|_2}{\|g-g^\nu\|_2\mathbb{P}(\nu<\infty)^{\frac{1}{p}-\frac{1}{2}}}\Bigg(\frac{q_1}{p_1}\int_0^\infty t^{\frac{q_1}{p_1}-1}\Big(\mu_t(\chi_{\{\nu<\infty\}})\Big)^{q_1}dt\Bigg)^{\frac{1}{q_1}}
\\&=&\frac{C}{\mathbb{P}(\nu<\infty)^{\frac{1}{p}-\frac{1}{2}}}\Bigg(\frac{q_1}{p_1}\int_0^\infty t^{\frac{q_1}{p_1}-1}\chi_{[0,\mathbb{P}(\nu<\infty))}^{q_1}(t)dt\Bigg)^{\frac{1}{q_1}}
\\&=&\frac{C}{P(\nu<\infty)^{\frac{1}{p}-\frac{1}{2}}}\mathbb{P}(\nu<\infty)^{\frac{1}{p_1}}=C.
\end{eqnarray*} Set $h_0=h/C,$  then $\|h_0\|_{H_{p,q}^s}\leq1$.
Consequently,
$\|\varphi\|\geq|\varphi(h_0)|=\mathbb{E}(h_0g)=\mathbb{E}(h_0(g-g^\nu))=C^{-1}\mathbb{P}(\nu<\infty)^{\frac{1}{2}-\frac{1}{p}}\|g-g^\nu\|_2$.
Taking the supremum over all stopping times, we have
$\|g\|_{BMO_2(\alpha)}\leq C\|\varphi\|$. The proof of the theorem
is complete.

Now we investigate the case $0<p\leq1, 1<q<\infty$.

\begin{thm}
The dual space of $H_{p,q}^s$ is $BMO_{2,q}(\alpha)$, $(0<p\leq 1,
1<q<\infty, \alpha=\frac{1}{p}-1)$.
\end{thm}

{Proof.} Let $g\in BMO_{2,q}(\alpha)\subset L_2$, define
$\varphi_g(f)=\mathbb{E}(fg)$, $(f\in L_2)$. Similarly to the proof
of Theorem 4.1, by H\"{o}lder's inequality we have
\begin{eqnarray*}
|\varphi_g(f)|&=&|\sum\limits_{k\in\mathbb{Z}}\mu_k\mathbb{E}(a^kg)|=|\sum\limits_{k\in\mathbb{Z}}\mu_k\mathbb{E}(a^k(g-g^{\nu_k}))|
\\&\leq&\sum\limits_{k\in\mathbb{Z}}|\mu_k|\mathbb{E}\big(|a^k(g-g^{\nu_k})|\big)\leq\sum\limits_{k\in\mathbb{Z}}|\mu_k|\|a^k\|_2\|g-g^{\nu_k}\|_2
\\&\leq&\sum\limits_{k\in\mathbb{Z}}|\mu_k|\mathbb{P}(\nu_k<\infty)^{\frac{1}{2}-\frac{1}{p}}\|g-g^{\nu_k}\|_2
=A\sum\limits_{k\in\mathbb{Z}}2^k\mathbb{P}(\nu_k<\infty)^{\frac{1}{2}}\|g-g^{\nu_k}\|_2.
\end{eqnarray*}
By the definition of $\|\cdot\|_{BMO_{2,q}(\alpha)}$ and Theorem
3.1, then
$$|\varphi_g(f)|\leq A\Big(\sum\limits_{k\in\mathbb{Z}}\big(2^k\mathbb{P}(\nu_k<\infty)^{\frac{1}{p}}\big)^q\Big)^{\frac{1}{q}}\|g\|_{BMO_{2,q}(\alpha)}\leq C\|f\|_{H_{p,q}^s}\|g\|_{BMO_{2,q}(\alpha)}.$$
Thus $\varphi_g$ can be uniquely extended to a continuous functional
on $H_{p,q}^s$.

Conversely, if $\varphi\in (H_{p,q}^s)^*$, we know that there exists
$g\in L_2$ such that $\varphi(f)=\mathbb{E}(fg)$, $(f\in L_2)$. Let
$\{\nu_k\}_{k\in\mathbb{Z}}$ be an arbitrary stopping time sequence
such that
$\big\{2^k\mathbb{P}(\nu_k<\infty)^{\frac{1}{p}}\big\}_{k\in\mathbb{Z}}\in
l_q$ and $N$ be an arbitrary nonnegative integer. Let
$$h_k=\frac{|g-g^{\nu_k}|
\textrm{sign}(g-g^{\nu_k})}{\|g-g^{\nu_k}\|_2}, \ \ \
f=\sum\limits_{k=-N}^N2^k\mathbb{P}(\nu_k<\infty)^{\frac{1}{2}}(h_k-h_k^{\nu_k}).$$
For an arbitrary integer $k_0$ which satisfies $-N\leq k_0\leq N$
(for $k_0\leq-N$, let $G=0$ and $H=f$; for $k_0>N$, let $H=0$ and
$G=f$), let $$f=G+H,$$ where
$G=\sum\limits_{k=-N}^{k_0-1}2^k\mathbb{P}(\nu_k<\infty)^{\frac{1}{2}}(h_k-h_k^{\nu_k})$
and
$H=\sum\limits_{k=k_0}^N2^k\mathbb{P}(\nu_k<\infty)^{\frac{1}{2}}(h_k-h_k^{\nu_k})$.
Obviously $\|h_k\|_{2}=1$, and
$\|G\|_{2}\leq2\sum\limits_{k=-N}^{k_0-1}2^k\mathbb{P}(\nu_k<\infty)^{\frac{1}{2}}$.
By the sublinearity of the operator $s$, we have $s(f)\leq
s(G)+s(H)$. Let $\varepsilon=\frac{p}{2}$, then
$0<\varepsilon<\min(1,\frac{q}{p})$. We obtain
\begin{eqnarray*}
2^{k_0p}\mathbb{P}(s(G)>2^{k_0})^\varepsilon&\leq&2^{k_0p}\Big(\frac{1}{2^{2k_0}}\|s(G)\|_{2}^{2}\Big)^\varepsilon\leq
C\cdot2^{k_0(p-2\varepsilon)}\|G\|_{2}^{2\varepsilon}
\\&\leq&C\Big(\sum\limits_{k=-N}^{k_0-1}2^k\mathbb{P}(\nu_k<\infty)^{\frac{1}{2}}\Big)^{2\varepsilon}
=C\Big(\sum\limits_{k=-N}^{k_0-1}2^k\mathbb{P}(\nu_k<\infty)^{\frac{\varepsilon}{p}}\Big)^{p}
\\&\leq&C\sum\limits_{k=-N}^{k_0-1}\big(2^k\mathbb{P}(\nu_k<\infty)^{\frac{\varepsilon}{p}}\big)^{p}
\leq
C\sum\limits_{k=-\infty}^{k_0-1}\big(2^k\mathbb{P}(\nu_k<\infty)^{\frac{\varepsilon}{p}}\big)^{p}.
\end{eqnarray*}
On the other hand,
$$\{s(H)>0\}\subset\bigcup\limits_{k=k_0}^N\{\nu_k<\infty\}.$$ Then
for each $0<\varepsilon<1$, we have
\begin{eqnarray*}
2^{k_0\varepsilon p}\mathbb{P}(s(H)>2^{k_0})&\leq&2^{k_0\varepsilon
p}\mathbb{P}(s(H)>0)\leq2^{k_0\varepsilon
p}\sum\limits_{k=k_0}^N\mathbb{P}(\nu_k<\infty)
\\&=&2^{k_0\varepsilon p}\sum\limits_{k=k_0}^N2^{k\varepsilon p}\mathbb{P}(\nu_k<\infty)2^{-k\varepsilon p}\leq\sum\limits_{k=k_0}^N2^{k\varepsilon p}\mathbb{P}(\nu_k<\infty)
\\&=&\sum\limits_{k=k_0}^N\big(2^{k\varepsilon}\mathbb{P}(\nu_k<\infty)^{\frac{1}{p}}\big)^p\leq\sum\limits_{k=k_0}^\infty\big(2^{k\varepsilon}\mathbb{P}(\nu_k<\infty)^{\frac{1}{p}}\big)^p.
\end{eqnarray*}
By Lemma 2.4, we have $s(f)\in L_{p,q}$ and $\|s(f)\|_{p,q}\leq
C\|\{2^k\mathbb{P}(\nu_k<\infty)^{\frac{1}{p}}\}_{k\in\mathbb{Z}}\|_{l_q}$.
Thus $f\in H_{p,q}^s$ and $$\|f\|_{H_{p,q}^s}\leq
C\Big(\sum\limits_{k\in\mathbb{Z}}\big(2^k\mathbb{P}(\nu_k<\infty)^{\frac{1}{p}}\big)^q\Big)^{\frac{1}{q}}.$$
Therefore,
\begin{eqnarray*}
\sum\limits_{k=-N}^N2^k\mathbb{P}(\nu_k<\infty)^{\frac{1}{2}}\|g-g^{\nu_k}\|_2&=&\sum\limits_{k=-N}^N2^k\mathbb{P}(\nu_k<\infty)^{\frac{1}{2}}\mathbb{E}(h_k(g-g^{\nu_k}))
\\&=&\sum\limits_{k=-N}^N2^k\mathbb{P}(\nu_k<\infty)^{\frac{1}{2}}\mathbb{E}((h_k-h_k^{\nu_k})g)
\\&=&\mathbb{E}(fg)=\varphi(f)\leq\|f\|_{H_{p,q}^s}\|\varphi\|
\\&\leq&C\Big(\sum\limits_{k\in\mathbb{Z}}\big(2^k\mathbb{P}(\nu_k<\infty)^{\frac{1}{p}}\big)^q\Big)^{\frac{1}{q}}\|\varphi\|.
\end{eqnarray*}
Thus we obtain
$$\frac{\sum\limits_{k=-N}^N2^k\mathbb{P}(\nu_k<\infty)^{\frac{1}{2}}\|g-g^{\nu_k}\|_2}{\Big(\sum\limits_{k\in\mathbb{Z}}\big(2^k\mathbb{P}(\nu_k<\infty)^{\frac{1}{p}}\big)^q\Big)^{\frac{1}{q}}}\leq
C\|\varphi\|.$$ Taking over all $N\in\mathbb{N}$ and the supremum
over all of such stopping time sequences satisfying
$\big\{2^k\mathbb{P}(\nu_k<\infty)^{\frac{1}{p}}\big\}_{k\in\mathbb{Z}}\in
l_q$, we get $\|g\|_{BMO_{2,q}(\alpha)}\leq C\|\varphi\|$. The proof
is complete.


\section{The generalized John-Nirenberg theorem}


\quad \quad In this section, we prove the generalized John-Nirenberg
theorem by duality when the stochastic basis
$\{\mathcal{F}_n\}_{n\geq0}$ is regular. Some of the dual results
are of independent interest. In order to do this, we need the
following lemma and we refer to \cite{W2} for these facts.

\begin{lem}
If the stochastic basis $\{\mathcal{F}_n\}_{n\geq0}$ is regular,
then the martingale Hardy-Lorentz spaces $H_{p,q}^*, H_{p,q}^S,
H_{p,q}^s, Q_{p,q}$ and $D_{p,q}$ are all equivalent for
$0<p<\infty, 0<q\leq\infty$, and $H_{p,q}^*, H_{p,q}^S, H_{p,q}^s,
Q_{p,q}, D_{p,q}$ and $L_{p,q}$ are all equivalent for $1<p<\infty,
0<q\leq\infty$.
\end{lem}

\begin{thm}
If the stochastic basis $\{\mathcal{F}_n\}_{n\geq0}$ is regular,
then
$$(H_{p,q}^s)^*=BMO_{r,q}(\alpha), \ \ \ (0<p\leq 1, 1<q,r<\infty, \alpha=\frac{1}{p}-1).$$
\end{thm}

{ Proof.} Let $g\in BMO_{r,q}(\alpha)\subset L_r$ and $r'$ be the
conjugate number of $r$, then $1<r'<\infty$. Define
$\varphi_g(f)=\mathbb{E}(fg)$, $f\in L_{r'}$. Note that
$L_{r'}=H^s_{r'} \subset H^s_{p,q} $. By Theorem 3.1 there exists a
sequence $(a^k)_{k\in\mathbb{Z}}$ of $(1,p,\infty)$-atoms and a
sequence of real numbers $(\mu_k)_{k\in\mathbb{Z}}$ satisfying
$\mu_k=A\cdot2^k\mathbb{P}(\nu_k<\infty)^{\frac{1}{p}}$ (where $A$
is a positive constant and $(\nu_k)_{k\in\mathbb{Z}}$ is the
corresponding stopping time sequence) such that
$f=\sum\limits_{k\in\mathbb{Z}}\mu_ka^k$ and
$\|(\mu_k)_{k\in\mathbb{Z}}\|_{l_q}\leq C\|f\|_{H_{p,q}^s}$. By
H\"{o}lder's inequality we can obtain
\begin{eqnarray*}
|\varphi_g(f)|&=&|\sum\limits_{k\in\mathbb{Z}}\mu_k\mathbb{E}(a^kg)|=|\sum\limits_{k\in\mathbb{Z}}\mu_k\mathbb{E}(a^k(g-g^{\nu_k}))|\leq\sum\limits_{k\in\mathbb{Z}}|\mu_k|\mathbb{E}\big(|a^k(g-g^{\nu_k})|\big)
\\&\leq&\sum\limits_{k\in\mathbb{Z}}|\mu_k|\|a^k\|_{r'}\|g-g^{\nu_k}\|_r\leq C\sum\limits_{k\in\mathbb{Z}}|\mu_k|\|s(a^k)\|_{r'}\|g-g^{\nu_k}\|_r
\\&\leq&C\sum\limits_{k\in\mathbb{Z}}|\mu_k|\mathbb{P}(\nu_k<\infty)^{\frac{1}{r'}-\frac{1}{p}}\|g-g^{\nu_k}\|_r=C\cdot A\sum\limits_{k\in\mathbb{Z}}2^k\mathbb{P}(\nu_k<\infty)^{1-\frac{1}{r}}\|g-g^{\nu_k}\|_r.
\end{eqnarray*}
By the definition of $\|\cdot\|_{BMO_{r,q}(\alpha)}$, we obtain
$$|\varphi_g(f)|\leq C\cdot A\Big(\sum\limits_{k\in\mathbb{Z}}
\big(2^k\mathbb{P}(\nu_k<\infty)^{\frac{1}{p}}\big)^q\Big)^{\frac{1}{q}}\|g\|_{BMO_{r,q}(\alpha)}\leq
C\|f\|_{H_{p,q}^s}\|g\|_{BMO_{r,q}(\alpha)}.$$ Thus $\varphi_g$ can
be extended to a continuous functional on $H_{p,q}^s$.

Conversely, if $\varphi\in(H_{p,q}^s)^*$. By the regularity of the
stochastic basis $\{\mathcal{F}_n\}_{n\geq0}$, we have $L_{r'}=
H_{r',r'}^s\subset H_{p,q}^s$, then
$(H_{p,q}^s)^*\subset(L_{r'})^*=L_r$. Thus there exists $g\in L_r$
such that $\varphi(f)=\varphi_g(f)=\mathbb{E}(fg)$, $(f\in L_{r'})$.

Let $\{\nu_k\}_{k\in\mathbb{Z}}$ be an arbitrary stopping time
sequence such that
$\big\{2^k\mathbb{P}(\nu_k<\infty)^{\frac{1}{p}}\big\}_{k\in\mathbb{Z}}\in
l_q$ and $N$ be an arbitrary nonnegative integer. Let
$$h_k=\frac{|g-g^{\nu_k}|^{r-1}
\textrm{sign}(g-g^{\nu_k})}{\|g-g^{\nu_k}\|_r^{r-1}}, \ \ \
f=\sum\limits_{k=-N}^N2^k\mathbb{P}(\nu_k<\infty)^{\frac{1}{r'}}(h_k-h_k^{\nu_k}).$$
For an arbitrary integer $k_0$ which satisfies $-N\leq k_0\leq N$
(for $k_0\leq-N$, let $G=0$ and $H=f$; for $k_0>N$, let $H=0$ and
$G=f$), let $$f=G+H,$$ where
$G=\sum\limits_{k=-N}^{k_0-1}2^k\mathbb{P}(\nu_k<\infty)^{\frac{1}{r'}}(h_k-h_k^{\nu_k})$
and
$H=\sum\limits_{k=k_0}^N2^k\mathbb{P}(\nu_k<\infty)^{\frac{1}{r'}}(h_k-h_k^{\nu_k})$.
Obviously, $\|h_k\|_{r'}=1$ and
$\|G\|_{r'}\leq2\sum\limits_{k=-N}^{k_0-1}2^k\mathbb{P}(\nu_k<\infty)^{\frac{1}{r'}}$.
By the sublinearity of the operator $s$, we have $s(f)\leq
s(G)+s(H)$. Let $\varepsilon=\frac{p}{r'}$, then
$0<\varepsilon<\min(1,\frac{q}{p})$. By Lemma 5.1 we have
\begin{eqnarray*}
2^{k_0p}\mathbb{P}(s(G)>2^{k_0})^\varepsilon&\leq&2^{k_0p}\Big(\frac{1}{2^{k_0r'}}\|s(G)\|_{r'}^{r'}\Big)^\varepsilon\leq
C\cdot2^{k_0(p-r'\varepsilon)}\|G\|_{r'}^{r'\varepsilon}
\\&\leq&C\Big(\sum\limits_{k=-N}^{k_0-1}2^k\mathbb{P}(\nu_k<\infty)^{\frac{1}{r'}}\Big)^{r'\varepsilon}
=C\Big(\sum\limits_{k=-N}^{k_0-1}2^k\mathbb{P}(\nu_k<\infty)^{\frac{\varepsilon}{p}}\Big)^{p}
\\&\leq&C\sum\limits_{k=-N}^{k_0-1}\big(2^k\mathbb{P}(\nu_k<\infty)^{\frac{\varepsilon}{p}}\big)^{p}
\leq
C\sum\limits_{k=-\infty}^{k_0-1}\big(2^k\mathbb{P}(\nu_k<\infty)^{\frac{\varepsilon}{p}}\big)^{p}.
\end{eqnarray*}
On the other hand,
$\{s(H)>0\}\subset\bigcup\limits_{k=k_0}^N\{\nu_k<\infty\}$. Then
for each $0<\varepsilon<1$, we have
\begin{eqnarray*}
2^{k_0\varepsilon p}\mathbb{P}(s(H)>2^{k_0})&\leq&2^{k_0\varepsilon
p}\mathbb{P}(s(H)>0)\leq2^{k_0\varepsilon
p}\sum\limits_{k=k_0}^N\mathbb{P}(\nu_k<\infty)
\\&\leq&\sum\limits_{k=k_0}^N2^{k\varepsilon p}\mathbb{P}(\nu_k<\infty)=\sum\limits_{k=k_0}^N\big(2^{k\varepsilon}\mathbb{P}(\nu_k<\infty)^{\frac{1}{p}}\big)^p
\\&\leq&\sum\limits_{k=k_0}^\infty\big(2^{k\varepsilon}\mathbb{P}(\nu_k<\infty)^{\frac{1}{p}}\big)^p.
\end{eqnarray*}
By Lemma 2.4, we have $s(f)\in L_{p,q}$ and $\|s(f)\|_{p,q}\leq
C\|\{2^k\mathbb{P}(\nu_k<\infty)^{\frac{1}{p}}\}_{k\in\mathbb{Z}}\|_{l_q}$.
Thus $f\in H_{p,q}^s$ and $$\|f\|_{H_{p,q}^s}\leq
C\Big(\sum\limits_{k\in\mathbb{Z}}\big(2^k\mathbb{P}(\nu_k<\infty)^{\frac{1}{p}}\big)^q\Big)^{\frac{1}{q}}.$$
Consequently,
\begin{eqnarray*}
\sum\limits_{k=-N}^N2^k\mathbb{P}(\nu_k<\infty)^{1-\frac{1}{r}}\|g-g^{\nu_k}\|_r&=&\sum\limits_{k=-N}^N2^k\mathbb{P}(\nu_k<\infty)^{\frac{1}{r'}}\mathbb{E}(h_k(g-g^{\nu_k}))
\\&=&\sum\limits_{k=-N}^N2^k\mathbb{P}(\nu_k<\infty)^{\frac{1}{r'}}\mathbb{E}((h_k-h_k^{\nu_k})g)
\\&=&\mathbb{E}(fg)=\varphi(f)\leq\|f\|_{H_{p,q}^s}\|\varphi\|
\\&\leq&C\Big(\sum\limits_{k\in\mathbb{Z}}\big(2^k\mathbb{P}(\nu_k<\infty)^{\frac{1}{p}}\big)^q\Big)^{\frac{1}{q}}\|\varphi\|.
\end{eqnarray*}
Thus we obtain
$$\frac{\sum\limits_{k=-N}^N2^k\mathbb{P}(\nu_k<\infty)^{1-\frac{1}{r}}\|g-g^{\nu_k}\|_r}{\Big(\sum\limits_{k\in\mathbb{Z}}\big(2^k\mathbb{P}(\nu_k<\infty)^{\frac{1}{p}}\big)^q\Big)^{\frac{1}{q}}}\leq
C\|\varphi\|.$$ Taking $N\rightarrow\infty$ and the supremum over
all of such stopping time sequences such that
$\big\{2^k\mathbb{P}(\nu_k<\infty)^{\frac{1}{p}}\big\}_{k\in\mathbb{Z}}\in
l_q$, we get $\|g\|_{BMO_{r,q}(\alpha)}\leq C\|\varphi\|$. The proof
is complete.

\bigskip It should be mentioned that the proof method of Theorem 5.2 is
not available for $r=1.$ In this case, we need new insight. Let the
dual space of $D_{p,q}$ be $D_{p,q}^*$. Let us denote by
$(D_{p,q}^*)_1$ those elements $\varphi$ from $D_{p,q}^*$ for which
there exists $g\in L_1$ such that $\varphi(f)=\mathbb{E}(fg)$, $f\in
L_\infty$. Namely
$$(D_{p,q}^*)_1=\{\varphi \in D_{p,q}^*: \ \exists g\in L_1 \ \ s.t. \ \ \varphi(f)=\mathbb{E}(fg), \ \forall f\in L_\infty\}.$$

\begin{thm}
$(D_{p,q}^*)_1=BMO_1(\alpha)$, $(0<p,q\leq 1,
\alpha=\frac{1}{p}-1)$.
\end{thm}
{ Proof.} Let $g\in BMO_1(\alpha)\subset L_1$. Define
$\varphi_g(f)=\mathbb{E}(fg)$, $(f\in L_\infty)$. By Theorem 3.3,
there exists a sequence $(a^k)_{k\in\mathbb{Z}}$ of
$(3,p,\infty)$-atoms and a sequence of real numbers
$(\mu_k)_{k\in\mathbb{Z}}$ satisfying
$\mu_k=A\cdot2^k\mathbb{P}(\nu_k<\infty)^{\frac{1}{p}}$ (where $A$
is a positive constant and $(\nu_k)_{k\in\mathbb{Z}}$ is the
corresponding stopping time sequence) such that
$f=\sum\limits_{k\in\mathbb{Z}}\mu_ka^k$ and
$\|(\mu_k)_{k\in\mathbb{Z}}\|_{l_q}\leq C\|f\|_{D_{p,q}}$. By
H\"{o}lder's inequality we obtain
\begin{eqnarray*}
|\varphi_g(f)|&=&|\sum\limits_{k\in\mathbb{Z}}\mu_k\mathbb{E}(a^kg)|=|\sum\limits_{k\in\mathbb{Z}}\mu_k\mathbb{E}(a^k(g-g^{\nu_k}))|
\\&\leq&\sum\limits_{k\in\mathbb{Z}}|\mu_k|\mathbb{E}\big(|a^k(g-g^{\nu_k})|\big)\leq\sum\limits_{k\in\mathbb{Z}}|\mu_k|\|a^k\|_\infty\|g-g^{\nu_k}\|_1
\\&\leq&\sum\limits_{k\in\mathbb{Z}}|\mu_k|\|(a^k)^*\|_\infty\|g-g^{\nu_k}\|_1\leq\sum\limits_{k\in\mathbb{Z}}|\mu_k|\mathbb{P}(\nu_k<\infty)^{-\frac{1}{p}}\|g-g^{\nu_k}\|_1
\\&\leq&\sum\limits_{k\in\mathbb{Z}}|\mu_k|\|g\|_{BMO_1(\alpha)}.
\end{eqnarray*}
Since $0<q\leq1$, then
$$|\varphi_g(f)|\leq\Big(\sum\limits_{k\in\mathbb{Z}}|\mu_k|^q\Big)^{\frac{1}{q}}\|g\|_{BMO_1(\alpha)}
\leq C\|f\|_{D_{p,q}}\|g\|_{BMO_1(\alpha)}.$$ Then $\varphi_g$ can
be extended to a continuous functional on $D_{p,q}$, and
$\varphi_g\in (D_{p,q}^*)_1$.

To prove the converse, let $\varphi\in(D_{p,q}^*)_1$, then there
exists $g\in L_1$ such that $\varphi(f)=\mathbb{E}(fg)$, $(f\in
L_\infty)$. Let $h=$sign$(g-g^\nu)$,
$a=\frac{1}{2}\mathbb{P}(\nu<\infty)^{-\frac{1}{p}}(h-h^\nu)$, where
$\nu\in\mathcal{T}$ is an arbitrary stopping time. Then $a$ is a
$(3,p,\infty)$-atom.

 Let $\mu=2A\cdot\mathbb{P}(\nu<\infty)^{\frac{1}{p}}$, let $h_0=\mu a=A(h-h^\nu)$. Considering the atomic decomposition of $h_0$, by Theorem 3.2 we have $h_0\in D_{p,q}$ and $\|h_0\|_{D_{p,q}}\leq C|\mu|=2CA\cdot\mathbb{P}(\nu<\infty)^{\frac{1}{p}}$, then $\|h-h^\nu\|_{D_{p,q}}\leq 2C\cdot\mathbb{P}(\nu<\infty)^{\frac{1}{p}}$. Thus we have
\begin{eqnarray*}
\mathbb{P}(\nu<\infty)^{-\frac{1}{p}}\|g-g^\nu\|_1&=&\mathbb{P}(\nu<\infty)^{-\frac{1}{p}}\mathbb{E}(h(g-g^\nu))
=\mathbb{P}(\nu<\infty)^{-\frac{1}{p}}\mathbb{E}((h-h^\nu)g)
\\&=&\mathbb{P}(\nu<\infty)^{-\frac{1}{p}}\varphi(h-h^\nu)
\leq
\mathbb{P}(\nu<\infty)^{-\frac{1}{p}}\|h-h^\nu\|_{D_{p,q}}\|\varphi\|
\\&=&2C\|\varphi\|.
\end{eqnarray*}
Taking the supremum over all stopping times, then we obtain
$\|g\|_{BMO_1(\alpha)}\leq C\|\varphi\|$. The proof of the theorem
is complete.

Now we consider $(D_{p,q}^*)_1$, $(0<p\leq1,1<q<\infty)$. We have
the following theorem.

\begin{thm}
$(D_{p,q}^*)_1=BMO_{1,q}(\alpha)$, $(0<p\leq 1, 1<q<\infty,
\alpha=\frac{1}{p}-1)$.
\end{thm}
{ Proof.} Let $g\in BMO_{1,q}(\alpha)\subset L_1$, then
$$\|g\|_{BMO_{1,q}(\alpha)}=\sup\frac{\sum\limits_{k\in\mathbb{Z}}2^k \|g-g^{\nu_k}\|_1}{\Big(\sum\limits_{k\in\mathbb{Z}}\big(2^k\mathbb{P}(\nu_k<\infty)^{\frac{1}{p}}\big)^q\Big)^{\frac{1}{q}}}<\infty,$$
where the supremum is taken over all stopping time sequences
$\{\nu_k\}_{k\in\mathbb{Z}}\subset\mathcal{T}$ such that $\big\{2^k
\mathbb{P}(\nu_k<\infty)^{\frac{1}{p}}\big\}_{k\in\mathbb{Z}}\in
l_q$. Define $\varphi_g(f)=\mathbb{E}(fg)$, $(f\in L_\infty)$.
Similarly to the proof of Theorem 4.3, by H\"{o}lder's inequality we
can obtain
\begin{eqnarray*}
|\varphi_g(f)|&=&|\sum\limits_{k\in\mathbb{Z}}\mu_k\mathbb{E}(a^kg)|=|\sum\limits_{k\in\mathbb{Z}}\mu_k\mathbb{E}(a^k(g-g^{\nu_k}))|
\\&\leq&\sum\limits_{k\in\mathbb{Z}}|\mu_k|\mathbb{E}\big(|a^k(g-g^{\nu_k})|\big)\leq\sum\limits_{k\in\mathbb{Z}}|\mu_k|\|a^k\|_\infty\|g-g^{\nu_k}\|_1
\\&\leq&\sum\limits_{k\in\mathbb{Z}}|\mu_k|\mathbb{P}(\nu_k<\infty)^{-\frac{1}{p}}\|g-g^{\nu_k}\|_1=A\sum\limits_{k\in\mathbb{Z}}2^k\|g-g^{\nu_k}\|_1.
\end{eqnarray*}
By the definition of $\|\cdot\|_{BMO_{1,q}(\alpha)}$, we obtain
$$|\varphi_g(f)|\leq A\Big(\sum\limits_{k\in\mathbb{Z}}\big(2^k\mathbb{P}(\nu_k<\infty)^{\frac{1}{p}}\big)^q\Big)^{\frac{1}{q}}\|g\|_{BMO_{1,q}(\alpha)}\leq C\|f\|_{D_{p,q}}\|g\|_{BMO_{1,q}(\alpha)}.$$
Thus $\varphi_g$ can be extended to a continuous functional on
$D_{p,q}$. Moreover, $\varphi_g\in(D_{p,q}^*)_1$.

Conversely, if $\varphi\in (D_{p,q}^*)_1$, then there exists $g\in
L_1$ such that $\varphi(f)=\mathbb{E}(fg)$, $(f\in L_\infty)$. Let
$\{\nu_k\}_{k\in\mathbb{Z}}$ be an arbitrary stopping time sequence
such that
$\big\{2^k\mathbb{P}(\nu_k<\infty)^{\frac{1}{p}}\big\}_{k\in\mathbb{Z}}\in
l_q$. Let
$$h_k=\textrm{sign}(g-g^{\nu_k}), \ \ \ a^k=\frac{1}{2}(h_k-h_k^{\nu_k})\mathbb{P}(\nu_k<\infty)^{-\frac{1}{p}}.$$
then $a^k$ is a $(3,p,\infty)$-atom.

Let
$f^N=\sum\limits_{k=-N}^N2^{k+1}\mathbb{P}(\nu_k<\infty)^{\frac{1}{p}}a^k$,
where $N$ is an arbitrary nonnegative integer. By Theorem 3.3 we
have $f^N\in D_{p,q}$ and
$$\|f^N\|_{D_{p,q}}\leq C\Big(\sum\limits_{k=-N}^N\big(2^k\mathbb{P}(\nu_k<\infty)^{\frac{1}{p}}\big)^q\Big)^{\frac{1}{q}}\leq C\Big(\sum\limits_{k\in\mathbb{Z}}\big(2^k\mathbb{P}(\nu_k<\infty)^{\frac{1}{p}}\big)^q\Big)^{\frac{1}{q}}.$$
Consequently,
\begin{eqnarray*}
\sum\limits_{k=-N}^N2^k\|g-g^{\nu_k}\|_1&=&\sum\limits_{k=-N}^N2^k\mathbb{E}(h_k(g-g^{\nu_k}))=\sum\limits_{k=-N}^N2^k\mathbb{E}((h_k-h_k^{\nu_k})g)
\\&=&\mathbb{E}(f^Ng)=\varphi(f^N)\leq\|f^N\|_{D_{p,q}}\|\varphi\|
\\&\leq&C\Big(\sum\limits_{k\in\mathbb{Z}}\big(2^k\mathbb{P}(\nu_k<\infty)^{\frac{1}{p}}\big)^q\Big)^{\frac{1}{q}}\|\varphi\|.
\end{eqnarray*}
Thus we have
$$\frac{\sum\limits_{k=-N}^N2^k\|g-g^{\nu_k}\|_1}{\Big(\sum\limits_{k\in\mathbb{Z}}\big(2^k\mathbb{P}(\nu_k<\infty)^{\frac{1}{p}}\big)^q\Big)^{\frac{1}{q}}}\leq
C\|\varphi\|.$$ This shows $\|g\|_{BMO_{1,q}(\alpha)}\leq
C\|\varphi\|$. The proof is complete.

\begin{prop}
If the stochastic basis $\{\mathcal{F}_n\}_{n\geq0}$ is regular, for
$0<p\leq1, 0<q<\infty$, then $(D_{p,q}^*)_1=D_{p,q}^*$.
\end{prop}
{ Proof.} Since $0<p\leq1$, then by Lemma 5.1, $L_2$ can also be
embedded continuously in $D_{p,q}$. Then $D_{p,q}^*\subset
(L_2)^*=L_2$. Let $\varphi$ be an arbitrary element of $D_{p,q}^*$,
then there exists $g\in L_2\subset L_1$ such that
$\varphi=\varphi_g$. By the definition of $(D_{p,q}^*)_1$, we have
$\varphi\in(D_{p,q}^*)_1$, then $D_{p,q}^*\subset(D_{p,q}^*)_1$. And
the inclusion relation $(D_{p,q}^*)_1\subset D_{p,q}^*$ is evident.
Hence we obtain $$(D_{p,q}^*)_1=D_{p,q}^*, \ \ \ (0<p\leq1,
0<q<\infty).$$ The proof of the proposition is complete.

We now are a position to prove Theorem 1.2.

\noindent { Proof of Theorem 1.2.}  It follows from Theorem 4.2 and
Theorem 5.2 that $$BMO_{r,q}(\alpha)=BMO_{2,q}(\alpha), \quad
1<r<\infty.$$ For $r=1,$ combining Theorem 4.2, Lemma 5.1, Theorem 5.4  with
Proposition 5.5, we get
$$BMO_{1,q}(\alpha)=BMO_{2,q}(\alpha).$$


\section{Boundedness of fractional integrals on martingale Hardy-Lorentz spaces}


\quad \quad As we know, Chao and Ombe \cite{CO} introduced the
fractional integrals for dyadic martingales. Recently, Nakai and
Sadasue \cite{NS} extended the notion of fractional integrals to
more general martingales. Sadasue \cite{S1} proves the boundedness
of fractional integrals on martingale Hardy spaces for $0<p\leq1$.
We now extend the boundedness of fractional integrals to martingale
Hardy-Lorentz spaces. In this section, we suppose that every
$\sigma$-algebra $\mathcal{F}_n$ is generated by countable atoms,
where $B\in\mathcal{F}_n$ is called an atom, if any $A\subset B$
with $A\in\mathcal{F}_n$ satisfies $\mathbb{P}(A)<\mathbb{P}(B)$,
then $\mathbb{P}(A)=0$. Denote by $A(\mathcal{F}_n)$ the set of all
atoms in $\mathcal{F}_n$. Without loss of generality, we always
suppose that the constant in (1.3) satisfying $R\geq2$.

Now we give the definition of fractional integral as follows.

\begin{defi}
For $f=(f_n)_{n\geq0}\in\mathcal{M}$, $\alpha>0$, the fractional
integral $I_\alpha f=\big((I_\alpha f)_n\big)_{n\geq0}$ of $f$ is
defined by $$(I_\alpha f)_n=\sum\limits_{k=1}^nb_{k-1}^\alpha
d_kf.$$ where $b_k$ is an $\mathcal{F}_k$-measurable function such
that $\forall B\in A(\mathcal{F}_k),\forall \omega\in B,
b_k(\omega)=\mathbb{P}(B)$.
\end{defi}
In order to prove the boundedness of fractional integrals, we need
the following lemmas.

\begin{lem}
Let $\{\mathcal{F}_n\}_{n\geq0}$ be regular, $f\in \mathcal{M}$ and
$\alpha>0$. Let $R$ be the constant in (1.3). If there exists
$B\in\mathcal{F}$ such that $f^*\leq \chi_B$. Then there exists a
positive constant
$C_\alpha$ independent of
$f$ and $B$ such that $$(I_\alpha f)^*\leq C_\alpha
\mathbb{P}(B)^\alpha\chi_B.$$
\end{lem}

For the proof of Lemma 6.2, see \cite{S1}, Lemma 3.5.

In the next lemma, we regard $(3,p,\infty)$-atom $a$ as a martingale
by $a=(a_n)_{n\geq0}=\big(E_n(a)\big)_{n\geq0}$, so we can consider
the fractional integral $I_\alpha a=\big((I_\alpha
a)_n\big)_{n\geq0}$.

\begin{lem}
Let $\{\mathcal{F}_n\}_{n\geq0}$ be regular and $R$ be the constant
in (1.3). If $0<p_1<p_2<\infty, \alpha=\frac{1}{p_1}-\frac{1}{p_2},
0<q_2\leq\infty$, and $a$ is a $(3,p_1,\infty)$-atom as in
Definition 2.2. Then we have
$$\|I_\alpha a\|_{H_{p_2,q_2}^*}\leq C_\alpha,$$ where $C_\alpha$ is the same constant as in Lemma 6.2.
\end{lem}

{ Proof.} Let $\nu$ be the stopping time associated with $a$. Then
we have $a^*\leq
\mathbb{P}(\nu<\infty)^{-\frac{1}{p_1}}\chi_{\{\nu<\infty\}}$.
Therefore
$\big(\mathbb{P}(\nu<\infty)^{\frac{1}{p_1}}a\big)^*=\mathbb{P}(\nu<\infty)^{\frac{1}{p_1}}a^*\leq\chi_{\{\nu<\infty\}}$.
By Lemma 6.2 we can obtain
$\big(I_\alpha(\mathbb{P}(\nu<\infty)^{\frac{1}{p_1}}a)\big)^*\leq
C_\alpha \mathbb{P}(\nu<\infty)^\alpha\chi_{\{\nu<\infty\}}$. Then
$$(I_\alpha a)^*\leq C_\alpha \mathbb{P}(\nu<\infty)^\alpha \mathbb{P}(\nu<\infty)^{-\frac{1}{p_1}}\chi_{\{\nu<\infty\}}=C_\alpha \mathbb{P}(\nu<\infty)^{-\frac{1}{p_2}}\chi_{\{\nu<\infty\}}.$$
By Proposition 2.1, we have
$$\mu_t\big((I_\alpha a)^*\big)\leq\mu_t\big(C_\alpha \mathbb{P}(\nu<\infty)^{-\frac{1}{p_2}}\chi_{\{\nu<\infty\}}\big)=C_\alpha \mathbb{P}(\nu<\infty)^{-\frac{1}{p_2}}\chi_{[0, \mathbb{P}(\nu<\infty))}(t).$$
For $0<q_2<\infty$, then
\begin{eqnarray*}
\|I_\alpha a\|_{H_{p_2,q_2}^*}^{q_2}&=&\|(I_\alpha
a)^*\|_{p_2,q_2}^{q_2}=\frac{q_2}{p_2}\int_0^\infty
t^{\frac{q_2}{p_2}-1}\Big(\mu_t\big((I_\alpha a)^*\big)\Big)^{q_2}dt
\\&\leq&\frac{q_2}{p_2}\int_0^\infty t^{\frac{q_2}{p_2}-1}\Big(C_\alpha \mathbb{P}(\nu<\infty)^{-\frac{1}{p_2}}\chi_{[0,\mathbb{P}(\nu<\infty))}(t)\Big)^{q_2}dt
\\&=&\frac{q_2}{p_2}\int_0^{\mathbb{P}(\nu<\infty)}t^{\frac{q_2}{p_2}-1}C_\alpha^{q_2}\mathbb{P}(\nu<\infty)^{-\frac{q_2}{p_2}}dt
\\&=&C_\alpha^{q_2}.
\end{eqnarray*}
For $q_2=\infty$, then
\begin{eqnarray*}
\|I_\alpha a\|_{H_{p_2,\infty}^*}&=&\|(I_\alpha
a)^*\|_{p_2,\infty}=\sup\limits_{t>0}t^{\frac{1}{p_2}}\mu_t((I_\alpha)^*)
\\&\leq&\sup\limits_{t>0}C_\alpha t^{\frac{1}{p_2}}\mathbb{P}(\nu<\infty)^{-\frac{1}{p_2}}\chi_{[0,\mathbb{P}(\nu<\infty))}(t)
\\&=&C_\alpha.
\end{eqnarray*}
Therefore $\|I_\alpha a\|_{H_{p_2,q_2}^*}\leq C_\alpha$, where
$C_\alpha$ is the same constant as in Lemma 6.2. The proof of is
complete.

\begin{thm}
Let $(\Omega,\mathcal{F},P)$ be a complete and nonatomic probability
space, and $\{\mathcal{F}_n\}_{n\geq0}$ be a regular stochastic
basis, let $0<q_1\leq1, q_1\leq q_2, q_1\leq p_2, 0<p_1<p_2<\infty,
\alpha=\frac{1}{p_1}-\frac{1}{p_2}$, then there exists a constant
$C$ such that $$\|I_\alpha f\|_{H_{p_2,q_2}^*}\leq
C\|f\|_{H_{p_1,q_1}^*},$$ for all $f\in H_{p_1,q_1}^*$.
\end{thm}
{Proof.} For $f\in H_{p_1,q_1}^*$. Since
$\{\mathcal{F}_n\}_{n\geq0}$ is regular, by Theorem 3.3 and Lemma
5.1, there exists a sequence $(a^k)_{k\in \mathbb{Z}}$ of
$(3,p_1,\infty)$-atoms and and a real number sequence $(\mu_k)_{k\in
\mathbb{Z}}\in l_{q_1}$ such that$$f_n=\sum\limits_{k\in
\mathbb{Z}}\mu_ka_n^k, \ \ \ (n\in\mathbb{N}),$$ and
$$\|(\mu_k)_{k\in\mathbb{Z}}\|_{l_{q_1}}\leq
C\|f\|_{H_{p_1,q_1}^*}.$$ Then by Lemma 6.3, we have
\begin{eqnarray*}
\|I_\alpha f\|_{H_{{p_2,q_2}}^*}^{q_1}&=&\|(I_\alpha
f)^*\|_{p_2,q_2}^{q_1}=\|(I_\alpha(\sum\limits_{k\in
\mathbb{Z}}\mu_ka^k))^*\|_{p_2,q_2}^{q_1}
\\&\leq&\|\sum\limits_{k\in \mathbb{Z}}|\mu_k|(I_\alpha a^k)^*\|_{p_2,q_2}^{q_1}\leq C\|\sum\limits_{k\in \mathbb{Z}}|\mu_k|(I_\alpha a^k)^*\|_{p_2,q_1}^{q_1}
\\&\leq&C\sum\limits_{k\in \mathbb{Z}}|\mu_k|^{q_1}\|(I_\alpha a^k)^*\|_{p_2,q_1}^{q_1}\leq C\cdot C_\alpha^{q_1}\|(\mu_k)_{k\in\mathbb{Z}}\|_{l_{q_1}}^{q_1}
\\&\leq&C\|f\|_{H_{p_1,q_1}^*}^{q_1}.
\end{eqnarray*}
Thus we have $$\|I_\alpha f\|_{H_{p_2,q_2}^*}\leq
C\|f\|_{H_{p_1,q_1}^*}.$$ The proof of the theorem is complete.

\begin{remark}
In Theorem 6.4, if we consider the special case $p_1=q_1=p,
p_2=q_2=q$, then we obtain the boundedness of fractional integrals
on martingale Hardy spaces for $0<p\leq1$, Theorem 3.1 in \cite{S1}
due to Sadasue.
\end{remark}

\end{document}